\documentclass{article}
\usepackage[utf8]{inputenc}
\usepackage{amssymb,amsmath,latexsym}
\usepackage[english]{babel}
\usepackage{amscd,amsxtra}
\usepackage{amsthm}
\usepackage{verbatim}
\usepackage{graphicx}
\usepackage{color}
\usepackage{tikz}
\usepackage{float}
\usetikzlibrary{calc, intersections, through, backgrounds}
\usetikzlibrary{shapes.geometric}
\tikzstyle{point}=[ball color=white, circle, draw=black, inner sep=0.1cm]
\tikzstyle{red}=[ball color=red, circle, draw=black, inner sep=0.1cm]
\tikzstyle{blue}=[ball color=blue, circle, draw=black, inner sep=0.1cm]
\usepackage{caption}
\usepackage[colorinlistoftodos]{todonotes}
\usetikzlibrary{arrows.meta}
\usepackage{authblk}

\newtheorem{theorem}{Theorem}[section]

\newtheorem{definition}[theorem]{Definition}

\newtheorem{lemma}[theorem]{Lemma}

\title{Vertex spans of multilayered cycles}
\author[1]{Aljoša Šubašić \thanks{Corresponding author: aljsub@pmfst.hr, Rudjera Bo\v skovi\' ca 33, 21000 Split, Croatia}}
\author[1]{Tanja Vojkovi\'c}
\affil[1]{Faculty of Science, University of Split}

\begin{document}

\maketitle
\begin{abstract}
In this paper we are defining a special class of graphs called multilayered graphs and its subclass, multilayered cycles. For that subclass of graphs we are giving the values of all vertex spans (strong, direct, or Cartesian span). Surprisingly, our results reveal that, irrespective of the chosen movement rules, the span values only depend on the length of the individual cycles, not the number of layers, which holds significant implications.  
\end{abstract}

Keywords: safety distance, graph spans, strong span, direct span, Cartesian span, multilayered graph \\
AMS Subject Classification: 05C12, 	05C90

\section{Introduction and motivation}

\begin{figure}[ht]
    \centering
    \begin{tikzpicture}[scale=1]
\node (1) at (1,1) [point] {};
\node (2) at (4,1) [point] {};
\node (3) at (7,1) [point] {};
\node (4) at (2,2) [point] {};
\node (5) at (4,2) [point] {};
\node (6) at (6,2) [point] {};
\node (7) at (3,3) [point] {};
\node (8) at (4,3) [point] {};
\node (9) at (5,3) [point] {};
\node (10) at (1,4) [point] {};
\node (11) at (2,4) [point] {};
\node (12) at (3,4) [point] {};
\node (13) at (5,4) [point] {};
\node (14) at (6,4) [point] {};
\node (15) at (7,4) [point] {};
\node (16) at (3,5) [point] {};
\node (17) at (4,5) [point] {};
\node (18) at (5,5) [point] {};
\node (19) at (2,6) [point] {};
\node (20) at (4,6) [point] {};
\node (21) at (6,6) [point] {};
\node (22) at (1,7) [point] {};
\node (23) at (4,7) [point] {};
\node (24) at (7,7) [point] {};
\draw (1)--(2)--(3)--(15)--(24)--(23)--(22)--(10)--(1) (4)--(5)--(6)--(14)--(21)--(20)--(19)--(11)--(4) (7)--(8)--(9)--(13)--(18)--(17)--(16)--(12)--(7) (2)--(5)--(8) (10)--(11)--(12) (13)--(14)--(15) (17)--(20)--(23);
\draw (1)--(4)--(7) (3)--(6)--(9) (16)--(19)--(22) (18)--(21)--(24);
\end{tikzpicture} 
    \caption{"Twelwe men's morris" game is played on $MC_8^3$}
    \label{TMM}
\end{figure}
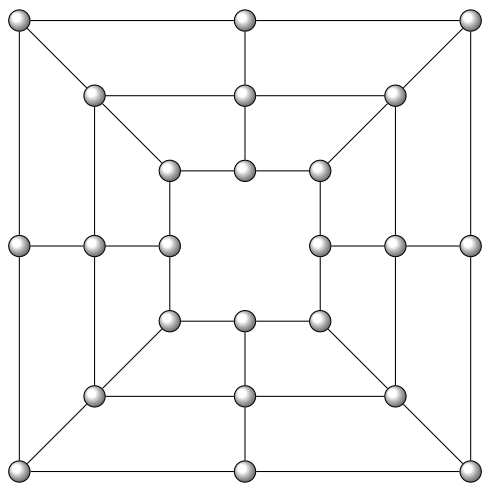

In response to the widespread adoption of social distancing measures during recent pandemic years, Banič and Taranenko used the concept of 'span' from topology, \cite{lelek}, and developed a graph theoretical measure known as the 'graph span', \cite{banic}. In its core it is the maximal safety distance two players can keep while moving through the vertices of some graph. Basic definitions introduce six types of graph spans, depending on weather the players have to visit all the vertices or all the edges of a graph (vertex and edge span) and what movement rules they follow (strong, direct and Cartesian span). In our previous research \cite{nas}, we explored the relationships between different types of vertex spans and determined span values for specific graph classes. Additionally, in \cite{nas2}, we conducted an analysis of edge spans and the minimum lengths of walks required to achieve these spans. In this paper we once again observe vertex spans, now for a special class of graphs - multilayered cycles, denoted by $MC_n^k$. We can imagine it as $k$ isomorphic cycles $C_n$, stacked on top of one another, with the corresponding vertices joined by "vertical "edges, forming a cylinder shape (Figure \ref{MLC}). Our results demonstrate that, regardless of the chosen movement rules, span values are solely dependent on the cycle's length, rather than the cylinder's height. This finding is particularly intriguing and has significant practical implications. Our motivation to observe this graph class stemmed from two distinct sources. Firstly, the need to determine safe occupancy limits in shopping malls, which are often designed in the multilayered cycle configuration, in response to social distancing measures. Secondly, our interest was piqued by graph-based games, such as 'The Cops and Robbers' games, \cite{igre}. There are different versions but the main idea is to have at least one "robber" and at least one "cop" moving through graph vertices, while cop is trying to "catch" a robber and the robber is trying to keep his distance from the cop. There are many versions of a game where multilayered cycles are a common playground, as observed from above it resembles a spider's web, for instance a game "The Spider and The Flies", developed in 1898. (Also "Web Chase" and "The Spider's Web"). It is also worth noting that a version of a well known game "Nine man's morris" called "Morabaraba" (or "Twelve men's morris"), which is played in South Africa as a sport, is played on a board that is exactly a multilayered cycle (Figure \ref{TMM}). 
In this paper, we provide a mathematical foundation for understanding span values in multilayered cycles. Future research could delve deeper into the mathematical analysis of the games inspired by these structures.
In Section \ref{prelim} basic definitions and preliminaries for our research are given and in Section \ref{results} we present our main results. Section \ref{conc} summerizes our results and presents some ideas for further work.

\section{Preliminaries and definitions}\label{prelim}

The term graph refers to a simple connected graph in the rest of the paper. We use standard graph theory notation, \cite{gross}.

\begin{definition}
    Let $G$ be a graph with $n$ vertices and $k\in \mathbb{N}, k \geq 2$. \textbf{Multilayered graph} $MG^k$, is a graph with $nk$ vertices, denoted by $(i,j)$, where $i \in V(G)$ and $j \in \mathbb{N}_k$, and in which $((a,b)(c,d)) \in E(MG^k)$ if one of the following holds:\\
    $a=c$ and $|d-b|=1$; or\\ 
    $b=d$ and $ac \in E(G)$.
\end{definition}

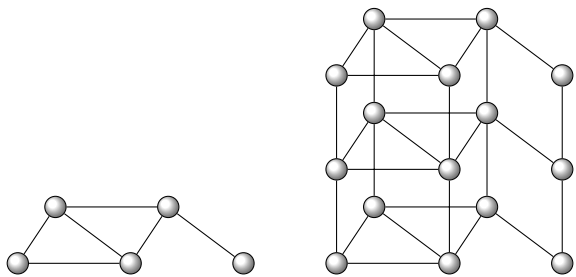
\begin{figure}[h!]
    \centering
    \begin{tikzpicture}[scale=0.5]
\node (1) at (1,0) [point] {};
\node (2) at (4,0) [point] {};
\node (3) at (7,0) [point] {};
\node (4) at (5,1.5) [point] {};
\node (5) at (2,1.5) [point] {};
\draw (1)--(2)--(4)--(5)--(1) (2)--(5) (4)--(3);
\end{tikzpicture} 
\hspace{20pt}
    \begin{tikzpicture}[scale=0.5]
\node (1) at (1,0) [point] {};
\node (2) at (4,0) [point] {};
\node (3) at (7,0) [point] {};
\node (4) at (5,1.5) [point] {};
\node (5) at (2,1.5) [point] {};
\node (6) at (1,2.5) [point] {};
\node (7) at (4,2.5) [point] {};
\node (8) at (7,2.5) [point] {};
\node (9) at (5,4) [point] {};
\node (10) at (2,4) [point] {};
\node (11) at (1,5) [point] {};
\node (12) at (4,5) [point] {};
\node (13) at (7,5) [point] {};
\node (14) at (5,6.5) [point] {};
\node (15) at (2,6.5) [point] {};
\draw (1)--(2)--(4)--(5)--(1) (2)--(5) (4)--(3) (6)--(7)--(9)--(10)--(6) (7)--(10) (9)--(8) (11)--(12)--(14)--(15)--(11) (12)--(15) (14)--(13) (1)--(6)--(11) (2)--(7)--(12) (3)--(8)--(13) (4)--(9)--(14) (5)--(10)--(15);
\end{tikzpicture} 
    \caption{Graphs $G$ and $MG^3$}
    \label{fig:MLG}
\end{figure}

The example of a multilayered graph is shown in Figure \ref{fig:MLG}. 

In this paper we will observe a special class of multilayered graphs, which are multilayered cycles. 

\begin{figure}[h!]
\centering
\begin{tikzpicture}
\node (1) at (1,0) [point] {0,1};
\node (2) at (4,0) [point] {1,1};
\node (3) at (6,1) [point] {2,1};
\node (4) at (5,1.5) [point] {3,1};
\node (5) at (2,1.5) [point] {4,1};
\node (6) at (0,0.5) [point] {5,1};
\node (7) at (1,2.5) [point] {0,2};
\node (8) at (4,2.5) [point] {1,2};
\node (9) at (6,3.5) [point] {2,2};
\node (10) at (5,4) [point] {3,2};
\node (11) at (2,4) [point] {4,2};
\node (12) at (0,3) [point] {5,2};
\node (13) at (1,5) [point] {0,3};
\node (14) at (4,5) [point] {1,3};
\node (15) at (6,6) [point] {2,3};
\node (16) at (5,6.5) [point] {3,3};
\node (17) at (2,6.5) [point] {4,3};
\node (18) at (0,5.5) [point] {5,3};
\draw (1) -- (2) -- (3) -- (4) -- (5) -- (6) -- (1) (7) -- (8) -- (9) -- (10) -- (11) -- (12) -- (7) (13) -- (14) -- (15) -- (16) -- (17) -- (18) -- (13) (1) -- (7) -- (13) (2) -- (8) -- (14) (3) -- (9) -- (15) (4) -- (10) -- (16) (5) -- (11) -- (17) (6) -- (12) -- (18);
\end{tikzpicture}
\caption{ Multilayered cycle $MC_6^3$}
\label{MLC}
\end{figure}
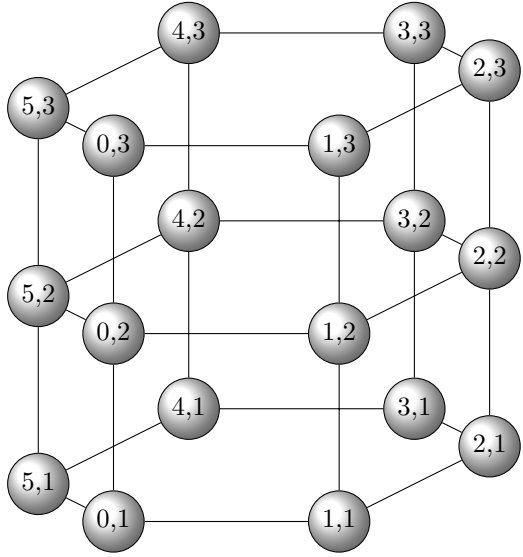

\begin{definition}
   Let $n,k\in \mathbb{N}, n \geq 3, k \geq 2$. \textbf{Multilayered cycle} $MC_n^k$, is a graph with $nk$ vertices, denoted by $(i,j)$, where $i \in \mathbb{Z}_n$ and $j \in \mathbb{N}_k$, and in which $((a,b)(c,d)) \in E(MC_n^k)$ if one of the following holds:\\
    $a=c$ and $|d-b|=1$; or\\ 
    $b=d$ and $|a-c| \in \{1,n-1\}$.
\end{definition} 
We will refer to the vertex $(i,j)$ as \textbf{the vertex $i$ in the layer $j$}. Also we will denote the layer of the vertex $(i,j)$ by $p_2(i,j):=j$, as layer is but a projection of a vertex to the second coordinate. The example of a multilayered cycle is given in Figure \ref{MLC}.

Note that a cube graph is the multilayered cycle $MC_4^2$.
To describe the movement of two players in s graph, we defined three types of functions that map $\mathbb{N}_l$, for some $l\in\mathbb{N}$, to the set of graph vertices, \cite{nas}. Such a function represents the movement of a player through graph vertices in $l$ steps.
These functions, and consequently, vertex spans, are defined  corresponding to three different movement rules that two players can apply in a graph so let us repeat those rules.
\begin{itemize}
    \item \textbf{Traditional movement rules}: Both players move independently of one another, one can stand still, while the other one moves, or they can both move at the same time;
    \item \textbf{Active movement rules}: Both players move to an adjacent vertex in each step;
    \item \textbf{Lazy movement rules}: In each step, exactly one of the players moves to an adjacent vertex while the other stands still.
\end{itemize}

The definitions of $l$-tracks, lazy $l$-tracks and opposite lazy $l$-tracks are given in \cite{nas}, but since they are of great importance for this paper, we will repeat them here.

\begin{definition}
	Let $G=(V,E)$ be a graph and $l\in\mathbb{N}$. We say that a surjective function $f_{l}:\mathbb{N}_{l}\longrightarrow V(G)$ is an\textbf{ $l$-track} on $G$ if $f(i)f(i+1)\in E(G)$ holds, for each $i\in\mathbb{N}_{l-1}$.
\end{definition}

 \begin{definition}\label{lazy}
    Let $G=(V,E)$ be a graph and $l\in\mathbb{N}$. We say that a surjective function $f_{l}:\mathbb{N}_{l}\longrightarrow V(G)$ is a\textbf{ lazy $l$-track} on $G$ if $f(i)f(i+1)\in E(G)$ or $f(i)=f(i+1)$ holds, for each $i\in\mathbb{N}_{l-1}$.
 \end{definition}

\begin{definition}
       Let $G$ be a graph, $f,g:\mathbb{N}_{l}\longrightarrow V(G)$ lazy $l$-tracks on $G$. We say that $f$ and $g$ are\textbf{ opposite lazy $l$-tracks} on $G$ if
       $$f(i)f(i+1)\in E(G) \iff g(i)=g(i+1),$$
       for all $i\in\mathbb{N}_{l-1}$.
\end{definition}

For reasons of simplifying our proofs, we will use the following terminology for lazy $l$-tracks in multilayered cycles, regarding the images of consequent steps. Let $MC_n^k$ be a multilayered cycle, $l\in\mathbb{N}$ and $f:\mathbb{N}_{l}\longrightarrow V(MC_n^k)$ a lazy $l$-track.
\begin{itemize}
    \item If $f(i)=(x,y)$ and $f(i+1)=(x,y)$ we will say that $f$ \textbf{stands still} in step $i$;
    \item If $f(i)=(x,y)$ and $f(i+1)=(x,y+1)$ we will say that $f$ \textbf{moves up} in step $i$;
    \item If $f(i)=(x,y)$ and $f(i+1)=(x,y-1)$ we will say that $f$ \textbf{moves down} in step $i$;
    \item If $f(i)=(x,y)$ and $f(i+1)=(x+_n1,y)$ we will say that $f$ \textbf{moves counter-clockwise} in step $i$;
    \item If $f(i)=(x,y)$ and $f(i+1)=(x-_n1,y)$ we will say that $f$ \textbf{moves clockwise} in step $i$.
\end{itemize}
Note that those are the only options for any lazy $l$-track on $MC_n^k$.

Next, we give the definition of the distance between two lazy $l$-tracks.

\begin{definition}
    Let $G$ be a graph, $l\in\mathbb{N}$ and $f,g$ two lazy $l$-tracks on $G$. We define
    $$m_G(f,g)=\min\{d(f(i),g(i):i\in\mathbb{N}_l\}$$
    to be the \textbf{distance between $f$ and $g$}.
\end{definition}

Analogously to an $l$-sweepable graph in \cite{nas2}, we define an $l$-trackable graph.

\begin{definition}
     Let $G$ be a graph and $l\in\mathbb{N}$. If at least one lazy $l$-track exists on $G$ we say that $G$ is an $l$-\textbf{trackable graph}.
\end{definition}

Lastly, we give definitions for different vertex spans, first described in \cite{banic}, and then redefined in \cite{nas}.

    Let $G$ be an $l$-trackable graph. We define
    
    $$M_{l}^{\boxtimes}:=\max\{m_{G}(f,g):f \text{ and } g \text{ are } \text{lazy }l  \text{-tracks on } G\}. $$
    
    $$M^{\times}_{l}:=\max\{m_{G}(f,g):f \text{ and } g \text{ are } l  \text{-tracks on } G\}. $$
 
    Let $G$ be a graph and $l\in\mathbb{N}$ such that at least one pair of opposite lazy $l$-tracks exists on $G$. We define
    $$M^{\square}_{l}:=\max\{m_{G}(f,g):f \text{ and } g \text{ are opposite lazy } l  \text{-tracks on } G\}.$$

    Let $G$ be a graph and let $S\subseteq\mathbb{N}$ be the set of all integers $l$ for which $G$ is an $l$-trackable graph. 
    We define the \textbf{strong vertex span} as the number
    
    $$\sigma^{\boxtimes}_V(G):=\max\{M_{l}^{\boxtimes}:l\in S\}.$$

    This number is the maximal safety distance that can be kept while two players visit all the vertices of a graph while following the traditional movement rules.

    We define the \textbf{direct vertex span} as the number
    
    $$\sigma^{\times}_V(G):=\max\{M^{\times}_{l}:l\in S\}.$$
    
    This number is the maximal safety distance that can be kept while two players visit all the edges of a graph with respect to the active movement rules. \\ 

    Let $G$ be a graph and let $C\subseteq\mathbb{N}$ be the set of all integers $l$ for which opposite lazy $l$-tracks exist on $G$. We define the \textbf{Cartesian vertex span} as the number
    
    $$\sigma^{\square}_V(G):=\max\{M^{\square}_{l}:l\in C\}.$$

    This number is the maximal safety distance that can be kept while two players visit all the edges of a graph with respect to the lazy movement rules.

\section{Results}\label{results}
We now proceed with the results for vertex spans values for multilayerd cycles. 

\begin{lemma} \label{Lema1}
Let graph $G=MC_n^k$, for some $n,k \in \mathbb{N}$. Also, let $f,g$ be two opposite lazy $l$-tracks on $G$, $l\in\mathbb{N}$. Then there exists $i \in \mathbb{N}_l$ such that $p_2(f(i))=p_2(g(i))$, i.e. $f(i)$ and $g(i)$ are in the same layer. 
\end{lemma}

\begin{proof}
Let $f,g$ be any two opposite lazy $l$-tracks on $G$. If $f(1)$ and $g(1)$ are in the same layer then the claim holds. Otherwise, let us assume that $f(1)$ is in layer $x$ and $g(1)$ is in layer $y \neq x$. Without any loss of generality, we can assume that $x<y$, so $x-y<0$. Since $f$ is surjective there exists $j \in \mathbb{N}_l$ such that $f(j)$ is in the layer $k$. For such $j$ it holds that $p_2(f(j))-p_2(g(j))=k-p_2(g(j)) \geq 0$. Since $f$ and $g$ are opposite, if $p_2(f(b))-p_2(g(b))=a$, for some $a,b \in \mathbb{N}$, then $p_2(f(b+1))-p_2(g(b+1)) \in \{a-1,a,a+1\}$. Now, given the facts that $p_2(f(1))-p_2(g(1))<0$, $p_2(f(j))-p_2(g(j)) \geq 0$ and that the difference between layers changes by at most 1 for consequent steps, we know that there must exist some $i \in \{2,...,j\}$ such that $p_2(f(i))-p_2(g(i))=0$ and therefore $p_2(f(i))=p_2(g(i))$.
\end{proof}
Example of one such movement is presented in Figure \ref{Slika 2}.

\begin{center}
\begin{figure}[ht]
    \centering
    \begin{tikzpicture}
\draw (0,1) -- (7.5,1) (0,2) -- (7.5,2) (0,3) -- (7.5,3) (0,1.5) -- (7.5,1.5) (0,2.5) -- (7.5,2.5);
\node (f1) at (1,1.5) [red] {};
\node [below=1pt]  at (1,1.5) {$f(1)$};
\node (f2) at (1.5,1.5) [red] {};
\node (f3) at (2,1) [red] {};
\node (f4) at (2.5,1.5) [red] {};
\node (f5) at (3,2) [red] {};
\node (f6) at (3.5,2) [red] {};
\node (f7) at (4,2) [] {};
\node [below=1pt]  at (4,2) {$f(i)$};
\node (f8) at (4.5,2) [red] {};
\node (f9) at (5,2.5) [red] {};
\node (f10) at (5.5,2.5) [red] {};
\node (f11) at (6,2.5) [red] {};
\node (f12) at (6.5,3) [red] {};
\node [below=1pt] at (6.5,3) {$f(j)$};
\node (g1) at (1,2.5) [blue] {};
\node [above=1pt] at (1,2.5) {$g(1)$};
\node (g2) at (1.5,3) [blue] {};
\node (g3) at (2,3) [blue] {};
\node (g4) at (2.5,3) [blue] {};
\node (g5) at (3,3) [blue] {};
\node (g6) at (3.5,2.5) [blue] {};
\node (g7) at (4,2) [] {};
\node [above=1pt]  at (4,2) {$g(i)$};
\node (g8) at (4.5,1.5) [blue] {};
\node (g9) at (5,1.5) [blue] {};
\node (g10) at (5.5,1) [blue] {};
\node (g11) at (6,1.5) [blue] {};
\node (g12) at (6.5,1.5) [blue] {};
\coordinate (x) at (4,2);
 \fill[blue] (x) + (0, 0.15) arc (90:270:0.15);
 \fill[red] (x) + (0, -0.15) arc (270:450:0.15);
\end{tikzpicture}  
    \caption{An example of movement through layers of two opposite lazy $l$-tracks}
    \label{Slika 2}
\end{figure}
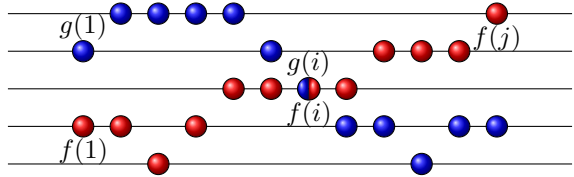
\end{center}

\begin{lemma} \label{Lema2}
$\sigma^{\square}_{V}(MC_n^k) \leq \left\lfloor\frac{n}{2}\right\rfloor$.
\end{lemma}

\begin{proof}
Since by Lemma \ref{Lema1}, for any two opposite lazy $l$-tracks $f$ and $g$, there exists an $i \in \mathbb{N}_l$ such that $p_2(f(i))=p_2(g(i))$, then for such $i$, $d(f(i),g(i)) \leq diam(C_n)=\left\lfloor\frac{n}{2}\right\rfloor$, hence $\sigma^{\square}_{V}(MC_n^k) \leq \left\lfloor\frac{n}{2}\right\rfloor$.
\end{proof}

\begin{theorem}\label{1}
 $\sigma^{\square}_{V}(MC_n^k)=\left\lfloor\frac{n}{2}\right\rfloor$.
\end{theorem}

\begin{proof}
First, we will construct two lazy $l$-tracks $f$ and $g$, where $l=2nk-1$, that are always at a distance at least $\left\lfloor\frac{n}{2}\right\rfloor$. We will start by defining $f(1)=(0,1)$ and $g(1)=(\left\lfloor\frac{n}{2}\right\rfloor,1)$, and then proceed in the following way: on odd steps $f$ moves and $g$ stands still, and on even steps $f$ stands still and $g$ moves. We will define all movement in four stages. 
\begin{itemize}
\item For the first $2k-2$ steps, on odd ones, $f$ moves up, and on even ones, $g$ moves up. So $f(2k-1)=(0,k)$ and $g(2k-1)=(\left\lfloor\frac{n}{2}\right\rfloor,k)$. 
\item On the next two steps, first $f$ moves counter-clockwise, and then $g$ moves counter-clockwise. So, $f(2k+1)=(1,k)$ and $g(2k+1)=(1+\left\lfloor\frac{n}{2}\right\rfloor,k)$. 
\item For the next $2k-2$ steps, on odd ones, $f$ moves down, and on even ones, $g$ moves down. So $f(4k-1)=(1,1)$ and $g(4k-1)=(1+\left\lfloor\frac{n}{2}\right\rfloor,1)$. 
\item On the next two steps, first $f$ moves counter-clockwise, and then $g$ moves counter-clockwise.
\end{itemize}
Now we repeat these four stages of movement until all vertices are visited. This kind of movement is presented in Figure \ref{fig:slika11}. It is easily seen that, defined this way, $f$ and $g$ are always at a distance of at least $\left\lfloor\frac{n}{2}\right\rfloor$. So, $\sigma^{\square}_{V}(MC_n^k) \geq \left\lfloor\frac{n}{2}\right\rfloor$. Combined with Lemma \ref{Lema2} we have $\sigma^{\square}_{V}(MC_n^k)=\left\lfloor\frac{n}{2}\right\rfloor$.
\begin{figure}[ht]
    \centering
\begin{tikzpicture}[scale=0.8] 
\node (1) at (1,0) [blue] {};
\node [below=1pt]  at (1,0) {$f(1)$};
\node (2) at (4,0) [blue] {};
\node (3) at (6,1) [blue] {};
\node (4) at (5,1.5) [red] {};
\node [below=1pt]  at (5,1.5) {$g(1)$};
\node (5) at (2,1.5) [red] {};
\node (6) at (0,0.5) [red] {};
\node (7) at (1,2.5) [blue] {};
\node (8) at (4,2.5) [blue] {};
\node (9) at (6,3.5) [blue] {};
\node (10) at (5,4) [red] {};
\node (11) at (2,4) [red] {};
\node (12) at (0,3) [red] {};
\node (13) at (1,5) [blue] {};
\node (14) at (4,5) [blue] {};
\node (15) at (6,6) [blue] {};
\node (16) at (5,6.5) [red] {};
\node (17) at (2,6.5) [red] {};
\node (18) at (0,5.5) [red] {};%
\draw [line width=1pt] (6) -- (1) -- (2) (3) -- (4) -- (5) (7) -- (8) -- (9) -- (10) -- (11) -- (12) -- (7) (14) -- (15) -- (16) (17) -- (18) -- (13);   
\draw[line width=2pt, red!, -to] (4) -- (10);
\draw[line width=2pt, red!, -to] (10) -- (16);
\draw[line width=2pt, red!, -to] (16) -- (17);
\draw[line width=2pt, red!, -to] (17) -- (11);
\draw[line width=2pt, red!, -to] (11) -- (5);
\draw[line width=2pt, red!, -to] (5) -- (6);
\draw[line width=2pt, red!, -to] (6) -- (12);
\draw[line width=2pt, red!, -to] (12) -- (18);
\draw[line width=2pt, blue!, -to] (1) -- (7);
\draw[line width=2pt, blue!, -to] (7) -- (13);
\draw[line width=2pt, blue!, -to] (13) -- (14);
\draw[line width=2pt, blue!, -to] (14) -- (8);
\draw[line width=2pt, blue!, -to] (8) -- (2);
\draw[line width=2pt, blue!, -to] (2) -- (3);
\draw[line width=2pt, blue!, -to] (3) -- (9);
\draw[line width=2pt, blue!, -to] (9) -- (15);
\end{tikzpicture} 
\captionof{figure}{ First half of the movement of $f$ and $g$ in Theorem \ref{1}} 
\label{fig:slika11}
\end{figure}
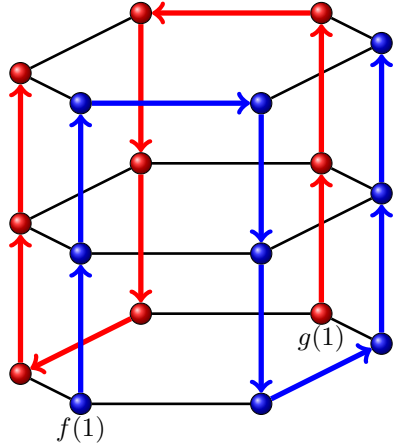
\end{proof}

\begin{lemma} \label{Lema3}
Let graph $G=MC_n^k$ for some $n,k \in \mathbb{N}$. Also, let $f,g$ be two lazy $l$-tracks on $G$, $l\in\mathbb{N}$. Then there exists $i \in \mathbb{N}_l$ such that $|p_2(f(i))-p_2(g(i))| \leq 1$, i.e. $f(i)$ and $g(i)$ are either in the same or in adjacent layers.
\end{lemma}

\begin{proof}
Much like in the proof of Lemma \ref{Lema1}, we can easily see that if, for some $b \in \mathbb{N}_l$, $|p_2(f(b))-p_2(g(b))|=a$, then $|p_2(f(b+1))-p_2(g(b+1))| \in \{a-2,a-1,a,a+1,a+2\}$. To put it in another words, in each step, the difference between layers of $f(b)$ and $g(b)$ can change by at most $2$. The same line of reasoning as in the proof of Lemma \ref{Lema1} leads us to conclusion that, for some $i \in \mathbb{N}_l$, $f(i)$ and $g(i)$ will be either in the same, or in the adjacent layers.
\end{proof}

Example of one such movement is presented in Figure \ref{fig:slika}.

\begin{center}
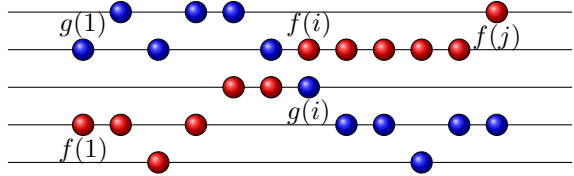
\begin{figure}[ht]
    \centering
    \begin{tikzpicture}
\draw (0,1) -- (7.5,1) (0,2) -- (7.5,2) (0,3) -- (7.5,3) (0,1.5) -- (7.5,1.5) (0,2.5) -- (7.5,2.5);
\node (f1) at (1,1.5) [red] {};
\node [below=1pt]  at (1,1.5) {$f(1)$};
\node (f2) at (1.5,1.5) [red] {};
\node (f3) at (2,1) [red] {};
\node (f4) at (2.5,1.5) [red] {};
\node (f5) at (3,2) [red] {};
\node (f6) at (3.5,2) [red] {};
\node (f7) at (4,2.5) [red] {};
\node [above=1pt]  at (4,2.5) {$f(i)$};
\node (f8) at (4.5,2.5) [red] {};
\node (f9) at (5,2.5) [red] {};
\node (f10) at (5.5,2.5) [red] {};
\node (f11) at (6,2.5) [red] {};
\node (f12) at (6.5,3) [red] {};
\node [below=1pt] at (6.5,3) {$f(j)$};
\node (g1) at (1,2.5) [blue] {};
\node [above=1pt] at (1,2.5) {$g(1)$};
\node (g2) at (1.5,3) [blue] {};
\node (g3) at (2,2.5) [blue] {};
\node (g4) at (2.5,3) [blue] {};
\node (g5) at (3,3) [blue] {};
\node (g6) at (3.5,2.5) [blue] {};
\node (g7) at (4,2) [blue] {};
\node [below=1pt]  at (4,2) {$g(i)$};
\node (g8) at (4.5,1.5) [blue] {};
\node (g9) at (5,1.5) [blue] {};
\node (g10) at (5.5,1) [blue] {};
\node (g11) at (6,1.5) [blue] {};
\node (g12) at (6.5,1.5) [blue] {};
\end{tikzpicture}  
    \caption{An example of movement through layers of two lazy $l$-tracks}
    \label{fig:slika}
\end{figure}
\end{center}

\begin{lemma} \label{Lema4}
$\sigma^{\boxtimes}_{V}(MC_n^k),\sigma^{\times}_{V}(MC_n^k) \leq \left\lfloor\frac{n}{2}\right\rfloor+1$.
\end{lemma}

\begin{proof}
By Lemma \ref{Lema3}, for any two lazy $l$-tracks $f$ and $g$, there exists an $i \in \mathbb{N}_l$ such that $f(i)$ and $g(i)$ are in the same or neighbouring layers. 
If $f(i)$ and $g(i)$ are in the same layer then $d(f(i),g(i)) \leq diam(C_n)=\left\lfloor\frac{n}{2}\right\rfloor$, hence both $\sigma^{\boxtimes}_{V}(MC_n^k)$ and $\sigma^{\times}_{V}(MC_n^k)$ are less than $\left\lfloor\frac{n}{2}\right\rfloor+1$. If, on the other hand, $f(i)$ and $g(i)$ are in the neighbouring layers then we can assume, without any loss of generality, that $f(i)$ is in one layer above $g(i)$, so $f(i)=(a,b)$ and $g(i)=(c,b-1)$ for some $a,b,c \in \mathbb{N}$. Now, $d(f(i),g(i))=d((a,b),(c,b))+1 \leq \left\lfloor\frac{n}{2}\right\rfloor+1$.
\end{proof}

\begin{theorem}
$\sigma^{\boxtimes}_{V}(MC_n^k)=\sigma^{\times}_{V}(MC_n^k)=\left\lfloor\frac{n}{2}\right\rfloor+1$.
\end{theorem}

\begin{proof}
We will construct two $l$-tracks $f$ and $g$ that will start on the distance $\left\lfloor\frac{n}{2}\right\rfloor+1$ and keep that distance at all times. First, let us describe the $l$-track $f$. It will start in vertex $(0,1)$, go through the whole layer clockwise, then go up one layer and go through the whole layer again in the same way. It will continue to do so until it reaches the topmost layer and goes through it as well. Lastly, it will go down one layer and once again go through it clockwise. For the graph 
$MC_6^4$, this movement is presented in Figure \ref{fig:movement of f}.
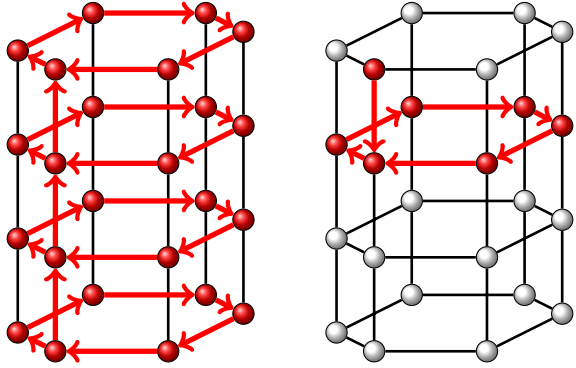
\begin{figure}[H]
    \centering
    \begin{tikzpicture}[scale=0.5]
\node (1) at (1,0) [red] {};
\node (2) at (4,0) [red] {};
\node (3) at (6,1) [red] {};
\node (4) at (5,1.5) [red] {};
\node (5) at (2,1.5) [red] {};
\node (6) at (0,0.5) [red] {};
\node (7) at (1,2.5) [red] {};
\node (8) at (4,2.5) [red] {};
\node (9) at (6,3.5) [red] {};
\node (10) at (5,4) [red] {};
\node (11) at (2,4) [red] {};
\node (12) at (0,3) [red] {};
\node (13) at (1,5) [red] {};
\node (14) at (4,5) [red] {};
\node (15) at (6,6) [red] {};
\node (16) at (5,6.5) [red] {};
\node (17) at (2,6.5) [red] {};
\node (18) at (0,5.5) [red] {};
\node (19) at (1,7.5) [red] {};
\node (20) at (4,7.5) [red] {};
\node (21) at (6,8.5) [red] {};
\node (22) at (5,9) [red] {};
\node (23) at (2,9) [red] {};
\node (24) at (0,8) [red] {};%
\draw [line width=1pt] (2)--(8)--(14)--(20) (3)--(9)--(15)--(21) (4)--(10)--(16)--(22) (5)--(11)--(17)--(23) (6)--(12)--(18)--(24);   
\draw[line width=2pt, red!, -to] (1) -- (6);
\draw[line width=2pt, red!, -to] (6) -- (5);
\draw[line width=2pt, red!, -to] (5) -- (4);
\draw[line width=2pt, red!, -to] (4) -- (3);
\draw[line width=2pt, red!, -to] (3) -- (2);
\draw[line width=2pt, red!, -to] (2) -- (1);
\draw[line width=2pt, red!, -to] (1) -- (7);
\draw[line width=2pt, red!, -to] (7) -- (12);
\draw[line width=2pt, red!, -to] (12) -- (11);
\draw[line width=2pt, red!, -to] (11) -- (10);
\draw[line width=2pt, red!, -to] (10) -- (9);
\draw[line width=2pt, red!, -to] (9) -- (8);
\draw[line width=2pt, red!, -to] (8) -- (7);
\draw[line width=2pt, red!, -to] (7) -- (13);
\draw[line width=2pt, red!, -to] (13) -- (18);
\draw[line width=2pt, red!, -to] (18) -- (17);
\draw[line width=2pt, red!, -to] (17) -- (16);
\draw[line width=2pt, red!, -to] (16) -- (15);
\draw[line width=2pt, red!, -to] (15) -- (14);
\draw[line width=2pt, red!, -to] (14) -- (13);
\draw[line width=2pt, red!, -to] (13) -- (19);
\draw[line width=2pt, red!, -to] (19) -- (24);
\draw[line width=2pt, red!, -to] (24) -- (23);
\draw[line width=2pt, red!, -to] (23) -- (22);
\draw[line width=2pt, red!, -to] (22) -- (21);
\draw[line width=2pt, red!, -to] (21) -- (20);
\draw[line width=2pt, red!, -to] (20) -- (19);

\end{tikzpicture} 
\hspace{20pt}
\begin{tikzpicture}[scale=0.5] 
\node (1) at (1,0) [point] {};
\node (2) at (4,0) [point] {};
\node (3) at (6,1) [point] {};
\node (4) at (5,1.5) [point] {};
\node (5) at (2,1.5) [point] {};
\node (6) at (0,0.5) [point] {};
\node (7) at (1,2.5) [point] {};
\node (8) at (4,2.5) [point] {};
\node (9) at (6,3.5) [point] {};
\node (10) at (5,4) [point] {};
\node (11) at (2,4) [point] {};
\node (12) at (0,3) [point] {};
\node (13) at (1,5) [red] {};
\node (14) at (4,5) [red] {};
\node (15) at (6,6) [red] {};
\node (16) at (5,6.5) [red] {};
\node (17) at (2,6.5) [red] {};
\node (18) at (0,5.5) [red] {};
\node (19) at (1,7.5) [red] {};
\node (20) at (4,7.5) [point] {};
\node (21) at (6,8.5) [point] {};
\node (22) at (5,9) [point] {};
\node (23) at (2,9) [point] {};
\node (24) at (0,8) [point] {};%
\draw [line width=1pt] (2)--(8)--(14)--(20) (3)--(9)--(15)--(21) (4)--(10)--(16)--(22) (5)--(11)--(17)--(23) (6)--(12)--(18)--(24) (1)--(2)--(3)--(4)--(5)--(6)--(1) (7)--(8)--(9)--(10)--(11)--(12)--(7) (19)--(20)--(21)--(22)--(23)--(24)--(19) (1)--(7)--(13);   
\draw[line width=2pt, red!, -to] (13) -- (18);
\draw[line width=2pt, red!, -to] (18) -- (17);
\draw[line width=2pt, red!, -to] (17) -- (16);
\draw[line width=2pt, red!, -to] (16) -- (15);
\draw[line width=2pt, red!, -to] (15) -- (14);
\draw[line width=2pt, red!, -to] (14) -- (13);
\draw[line width=2pt, red!, -to] (19) -- (13);
\end{tikzpicture} 
    \caption{The first and second part of movement of $f$ in graph $MC_6^4$}
    \label{fig:movement of f}
\end{figure}
It is easily seen that this way $f$ visits all the vertices. Now we will describe the movement of $g$ depending on the movement of $l$-track $f$. $l$-track $g$ will start its movement in vertex $(\left\lfloor\frac{n}{2}\right\rfloor,2)$, so $d(f(1),g(1))=\left\lfloor\frac{n}{2}\right\rfloor+1$.
Whenever $l$-track $f$ moves clockwise, $g$ also moves clockwise thus maintaining the same distance as well as visiting its whole layer while $f$ is visiting its own. The first time that $f$ moves up $g$ will move down, and afterwards whenever $f$ changes layers $g$ will move up. This way $g$ will get to visit all the layers, and at each one go through all of its vertices, while maintaining the same distance at all times. For the graph $MC_6^4$, the movement of $g$ is presented in Figure \ref{fig:movement of g}.
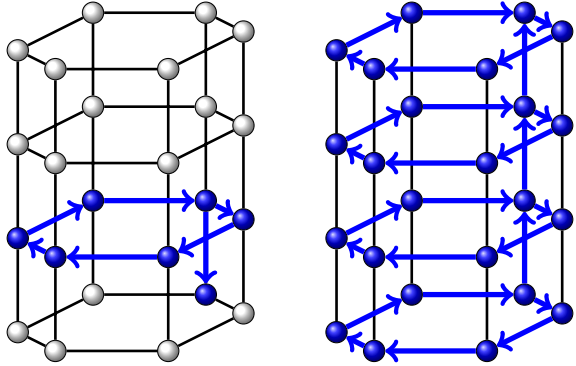
\begin{figure}[H]
    \centering
    \begin{tikzpicture}[scale=0.5]
\node (1) at (1,0) [point] {};
\node (2) at (4,0) [point] {};
\node (3) at (6,1) [point] {};
\node (4) at (5,1.5) [blue] {};
\node (5) at (2,1.5) [point] {};
\node (6) at (0,0.5) [point] {};
\node (7) at (1,2.5) [blue] {};
\node (8) at (4,2.5) [blue] {};
\node (9) at (6,3.5) [blue] {};
\node (10) at (5,4) [blue] {};
\node (11) at (2,4) [blue] {};
\node (12) at (0,3) [blue] {};
\node (13) at (1,5) [point] {};
\node (14) at (4,5) [point] {};
\node (15) at (6,6) [point] {};
\node (16) at (5,6.5) [point] {};
\node (17) at (2,6.5) [point] {};
\node (18) at (0,5.5) [point] {};
\node (19) at (1,7.5) [point] {};
\node (20) at (4,7.5) [point] {};
\node (21) at (6,8.5) [point] {};
\node (22) at (5,9) [point] {};
\node (23) at (2,9) [point] {};
\node (24) at (0,8) [point] {};%
\draw [line width=1pt] (1)--(7)--(13)--(19) (2)--(8)--(14)--(20) (3)--(9)--(15)--(21) (4)--(10)--(16)--(22) (5)--(11)--(17)--(23) (6)--(12)--(18)--(24) (1)--(2)--(3)--(4)--(5)--(6)--(1) (13)--(14)--(15)--(16)--(17)--(18)--(13) (19)--(20)--(21)--(22)--(23)--(24)--(19);   
\draw[line width=2pt, blue!, -to] (7) -- (12);
\draw[line width=2pt, blue!, -to] (12) -- (11);
\draw[line width=2pt, blue!, -to] (11) -- (10);
\draw[line width=2pt, blue!, -to] (10) -- (9);
\draw[line width=2pt, blue!, -to] (9) -- (8);
\draw[line width=2pt, blue!, -to] (8) -- (7);
\draw[line width=2pt, blue!, -to] (10) -- (4);
\end{tikzpicture} 
\hspace{20pt}
    \begin{tikzpicture}[scale=0.5]
\node (1) at (1,0) [blue] {};
\node (2) at (4,0) [blue] {};
\node (3) at (6,1) [blue] {};
\node (4) at (5,1.5) [blue] {};
\node (5) at (2,1.5) [blue] {};
\node (6) at (0,0.5) [blue] {};
\node (7) at (1,2.5) [blue] {};
\node (8) at (4,2.5) [blue] {};
\node (9) at (6,3.5) [blue] {};
\node (10) at (5,4) [blue] {};
\node (11) at (2,4) [blue] {};
\node (12) at (0,3) [blue] {};
\node (13) at (1,5) [blue] {};
\node (14) at (4,5) [blue] {};
\node (15) at (6,6) [blue] {};
\node (16) at (5,6.5) [blue] {};
\node (17) at (2,6.5) [blue] {};
\node (18) at (0,5.5) [blue] {};
\node (19) at (1,7.5) [blue] {};
\node (20) at (4,7.5) [blue] {};
\node (21) at (6,8.5) [blue] {};
\node (22) at (5,9) [blue] {};
\node (23) at (2,9) [blue] {};
\node (24) at (0,8) [blue] {};%
\draw [line width=1pt] (2)--(8)--(14)--(20) (3)--(9)--(15)--(21) (1)--(7)--(13)--(19) (5)--(11)--(17)--(23) (6)--(12)--(18)--(24);   
\draw[line width=2pt, blue!, -to] (1) -- (6);
\draw[line width=2pt, blue!, -to] (6) -- (5);
\draw[line width=2pt, blue!, -to] (5) -- (4);
\draw[line width=2pt, blue!, -to] (4) -- (3);
\draw[line width=2pt, blue!, -to] (3) -- (2);
\draw[line width=2pt, blue!, -to] (2) -- (1);
\draw[line width=2pt, blue!, -to] (4) -- (10);
\draw[line width=2pt, blue!, -to] (7) -- (12);
\draw[line width=2pt, blue!, -to] (12) -- (11);
\draw[line width=2pt, blue!, -to] (11) -- (10);
\draw[line width=2pt, blue!, -to] (10) -- (9);
\draw[line width=2pt, blue!, -to] (9) -- (8);
\draw[line width=2pt, blue!, -to] (8) -- (7);
\draw[line width=2pt, blue!, -to] (10) -- (16);
\draw[line width=2pt, blue!, -to] (13) -- (18);
\draw[line width=2pt, blue!, -to] (18) -- (17);
\draw[line width=2pt, blue!, -to] (17) -- (16);
\draw[line width=2pt, blue!, -to] (16) -- (15);
\draw[line width=2pt, blue!, -to] (15) -- (14);
\draw[line width=2pt, blue!, -to] (14) -- (13);
\draw[line width=2pt, blue!, -to] (16) -- (22);
\draw[line width=2pt, blue!, -to] (19) -- (24);
\draw[line width=2pt, blue!, -to] (24) -- (23);
\draw[line width=2pt, blue!, -to] (23) -- (22);
\draw[line width=2pt, blue!, -to] (22) -- (21);
\draw[line width=2pt, blue!, -to] (21) -- (20);
\draw[line width=2pt, blue!, -to] (20) -- (19);
\end{tikzpicture} 
    \caption{The first and second part of movement of $g$  in graph $MC_6^4$}
    \label{fig:movement of g}
\end{figure}
It is easily seen that, defined this way, $f$ and $g$ are always at a distance of at least $\left\lfloor\frac{n}{2}\right\rfloor+1$. So, both $\sigma^{\boxtimes}_{V}(MC_n^k)$ and $\sigma^{\times}_{V}(MC_n^k)$ are greater or equal to $\left\lfloor\frac{n}{2}\right\rfloor+1$. Combined with Lemma \ref{Lema4} we have $\sigma^{\boxtimes}_{V}(MC_n^k)=\sigma^{\times}_{V}(MC_n^k)=\left\lfloor\frac{n}{2}\right\rfloor+1$.
\end{proof}

\section{Summary and conclusion}\label{conc}

Let $MC_n^k$ be a multilayerd cycle. The summary of our results is given in Table \ref{summ}.
\renewcommand{\arraystretch}{1.5}
\begin{center}
\captionof{table}{Vertex span values for $MC_n^k$}\label{summ}
\begin{tabular}{ |c|c|c| } 
 \hline
 $\sigma^{\boxtimes}_{V}$ & $\sigma^{\times}_{V}$ & $\sigma^{\square}_{V}$ \\ 
 \hline
 $\left\lfloor\frac{n}{2}\right\rfloor+1$ & $\left\lfloor\frac{n}{2}\right\rfloor+1$ & $\left\lfloor\frac{n}{2}\right\rfloor$ \\ 

 \hline
\end{tabular}
\end{center}

In conclusion, we see that all vertex spans are depending only on the size of the cycle, and not on the number of layers.\\
Further work on multilayered graphs might include exploring the relation between $\sigma^{\square}_{V}(G),\sigma^{\times}_{V}(G),\sigma^{\boxtimes}_{V}(G)$ and $\sigma^{\square}_{V}(MG^k),\sigma^{\times}_{V}(MG^k),\sigma^{\boxtimes}_{V}(MG^k)$. We are also interested to find the edge spans for this graph class as well as all of the span values for the movement of more than two players. As stated in the introduction, it would be interesting to analyze the other aspects of games for which multilayered cycle is the playground.\\

\textbf{Author Contributions:} Conceptualization, A.Š. and T.V.; Investigation, A.Š. and T.V.; Data curation, A.Š. and T.V.; Writing—original draft, A.Š. and T.V.; Writing—review \& editing, A.Š. and T.V. All authors have read and agreed to the published version of the manuscript.\\

\textbf{Funding:} This research received no external funding.\\

\textbf{Data Availability Statement:} Not applicable.\\

\textbf{Acknowledgments:} Not applicable.\\

\textbf{Conflicts of Interest:} The authors declare no conflict of interest.

\end{document}